\documentclass[11pt]{article}
\usepackage[a4paper]{anysize}\marginsize{3.5cm}{3.5cm}{1.3cm}{2cm}
\pdfpagewidth=\paperwidth \pdfpageheight=\paperheight
\usepackage{amsfonts,amssymb,amsthm,amsmath,eucal}
\usepackage{pgf}
\usepackage{bbm,array}
\usepackage{tikz} 
\usepackage{float}
\restylefloat{table}
\usepackage{subfigure}
\usepackage{caption}
\usetikzlibrary{arrows}
\usepackage{multicol}
\usepackage{multirow}

\theoremstyle{plain}
\newtheorem{thm}{Theorem}[section]
\newtheorem{theorem}[thm]{Theorem}
\newtheorem*{theoremA}{Theorem A}
\newtheorem*{theoremB}{Theorem B}

\newtheorem{conjecture}[thm]{Conjecture}

\theoremstyle{definition}
\newtheorem{definition}[thm]{Definition}

\newtheorem{thevarthm}[thm]{\varthmname}

\newenvironment{varthm*}[1]{\trivlist\item[]{\bf #1.}\it}{\endtrivlist}

\newcommand\be{\begin{eqnarray*}}
\newcommand\ee{\end{eqnarray*}}

\newcommand\Q{\mathbb Q}
\newcommand\R{\mathbb R}

\newcommand\K{\mathbb K}

\renewcommand\P{\mathbb P}

\newcommand\newop[2]{\def#1{\mathop{\rm #2}\nolimits}}
\newop\mod{mod}
\newop\log{log}
\newop\ord{ord}
\newop\Gal{Gal}
\newop\SL{SL}
\newop\GL{GL}
\newop\Bl{Bl}
\newop\mult{mult}
\newop\mass{mass}
\newop\div{div}
\newop\codim{codim}
\newop\sing{sing}
\newop\vdim{vdim}
\newop\edim{edim}
\newop\Ass{Ass}
\newop\size{size}
\newop\reg{reg}
\newop\areg{areg}
\newop\asreg{asreg}
\newop\satdeg{satdeg}
\newop\supp{supp}
\newop\gin{gin}
\newop\ini{in}
\newop\vol{vol}
\newop\sat{sat}
\newop\length{length}
\newop\depth{depth}
\newop\characteristic{char}

\def\keywordname{{\bfseries Keywords}}%
\def\keywords#1{\par\addvspace\medskipamount{\rightskip=0pt plus1cm
\def\and{\ifhmode\unskip\nobreak\fi\ $\cdot$
}\noindent\keywordname\enspace\ignorespaces#1\par}}
\def\subclassname{{\bfseries Mathematics Subject Classification
(2010)}\enspace}
\def\subclass#1{\par\addvspace\medskipamount{\rightskip=0pt plus1cm
\def\and{\ifhmode\unskip\nobreak\fi\ $\cdot$
}\noindent\subclassname\ignorespaces#1\par}}

\begin{document}

\author{M. Janasz, M. Lampa-Baczy\'nska, G. Malara}
\title{New phenomena in the containment problem for simplicial arrangements}
\date{\today}
\maketitle
\thispagestyle{empty}

\begin{abstract}
In this note we consider two simplicial arrangements of lines and ideals $I$ of intersection points of these lines. There are $127$ intersection points in both cases and the numbers $t_i$ of points lying on exactly $i$ configuration lines (points of multiplicity $i$) coincide. We show that in one of these examples the containment $I^{(3)} \subseteq I^2$ holds, whereas it fails in the other. We also show that the containment fails for a subarrrangement of $21$ lines. The interest in the containment relation between $I^{(3)}$ and $I^2$ for ideals of points in $\P^2$ is motivated by a question posted by Huneke around $2000$. Configurations of points with $I^{(3)} \not\subseteq I^2$ are quite rare. Our example reveals two particular features: All points are defined over $\Q$ and all intersection points of lines are involved. In examples studied by now only points with multiplicity $i\geq 3$ were considered. The novelty of our arrangements lies in the geometry of the element in $I^{(3)}$ which witness the noncontainment in $I^2$. In all previous examples such an element was a product of linear forms. Now, in both cases there is an irreducible curve of higher degree involved.

 \subclass{14N20, 13A15, 13F20, 52C35 }
\end{abstract}


\section{Introduction}\label{intro}

 The following problem has attracted a lot of attention in the last two decades.

\vspace{0.2cm}
\begin{problem}
	
	\noindent Determine all pairs of positive integers $(m,r)$ such that the containment
	\begin{equation}
	\label{eq:cont}
		I^{(m)} \subset I^r	
	\end{equation}  
	holds for all homogeneous ideals $I\subseteq \mathbb{K}[x_0, \ldots ,x_{N}]$ in the ring  of polynomials over a field $\mathbb{K}$.
\end{problem}
\vspace{0.2cm}

In $2000$ Ein, Lazarsfeld and Smith \cite{ELS} in characteristic zero and Hochster and  Huneke \cite{HoHu} in positive characteristic discovered that the containment (\ref{eq:cont})  holds provided $m \geq Nr$. 

\begin{theorem}[Ein-Lazarsfeld-Smith, Hochster-Huneke]
	Let $I \subseteq \mathbb{K}[x_0, \dots, x_N]$ be a homogeneous ideal. Then there is 
	$$I^{(m)} \subset I^r$$
	for all $m \geq Nr$.
\end{theorem}

This ground-breaking result prompted a natural question about the optimality of the bound $m\geq Nr$. A number of examples suggested the following conjectural improvement (see \cite[Conjecture $8.4.2$]{primer},  or \cite[Conjecture $4.1.1$]{HaHu},  or \cite[Conjecture $1.1$]{BoHa})

\begin{conjecture} \label{symb}
	Let $I$ be a homogeneous ideal. Ten 
	$$I^{(m)} \subseteq I^r$$
	for all $m\geq Nr-(N-1)$.
\end{conjecture}

The first non-trivial case is $N=2$ and $r=2$.  Then there is always $I^{(4)}\subset I^2$ and it is very easy to give examples with $I^{(2)}\nsubseteq I^2$. Huneke asked around $2000$ if 
$$I^{(3)}\subseteq I^2$$
holds for all ideals defining points in $\P^2$.

This is not the case. The first non-containment example was announced in \cite{DST13} and soon after additional non-containment examples were discovered and described in \cite{CGMLLPS2015}, \cite{HarSec13}, \cite{SS}, \cite{MJ2015}, \cite{malara-szpond}, \cite{combinatorics}.

Such examples are quite rare and they all follow the same pattern, in particular 
they are related to line arrangements. More precisely, let $\mathcal{L}=\{L_1,\ldots,L_s\}$ be an arrangement of lines in $\P^2$ and let $\mathcal{P}=\{P_1,\ldots,P_t\}$ be the set of all points contained in at least $2$ lines from $\mathcal{L}$. Let $I$ be the ideal of those points which are contained in at least $3$ lines. By the Zariski-Nagata Theorem \cite[Theorem 3.14]{Eisenbud} the product 
$$f=l_1\cdot ... \cdot l_s\in I^{(3)}$$ 
and sometimes it happens that $f\not\in I^2$ (here $l_i$ is the equation of $L_i$). 

The novelty of our non-containment example is that whereas the ideal of points is determined by lines, it is not their product which sits in $I^{(3)}\setminus I^2.$ More precisely, our main results are the following

\begin{theoremA}  \label{confA}
	There exists an arrangement of $31$ lines which intersect in the total of $127$ points such that for the ideal $I$ of these $127$ points there is 
	\[I^{(3)}\nsubseteq I^2.\]
\end{theoremA}

\noindent
Moreover, there is an  element $f$ of degree $33$ in $I^{(3)}$, which is not contained in $I^2$ and which is a product of 
\begin{itemize}
	\item[$\bullet$] $21$ of arrangement lines and 
	\item[$\bullet$] an irreducible curve of degree $12$.
\end{itemize}


\begin{theoremB}
	There exists an arrangement of $21$ lines which intersect in the total of $115$ points such that for the ideal $I$ of these $115$ points there is 
	\[I^{(3)}\nsubseteq I^2.\]
\end{theoremB}
\noindent
Moreover, there is an  element $f$ of degree $31$ in $I^{(3)}$, which is not contained in $I^2$ and which is a product of 
\begin{itemize}
	\item[$\bullet$] all arrangement lines and 
	\item[$\bullet$] an irreducible curve of degree $10$.
\end{itemize}

A number of elementary but tedious calculations is omitted. Instead we provide a Singular script \cite{script} which provides easy verification of our claims.

\section{Preliminaries}
\label{sec:prelim}

In this section we define the basic object we are interested in and state the central conjecture in the field, which motivated our research here.

Let $I\subseteq \K[x_0,\ldots,x_N]$ be a homogeneous ideal in the ring of polynomials over a field $\K$.

\begin{definition} (Symbolic power)
	For $m\geq 1$, the $m$-th symbolic power of $I$ is the ideal 
	$$ I^{(m)} = \mathbb{K}[\mathbb{P}^N] \cap \left( \bigcap_{\mathfrak{p} \in \Ass(I)} (I^{m})_{\mathfrak{p}} \right), $$
	where the intersection is taken over all associated primes of $I$.
\end{definition}

Symbolic powers of ideals are of geometric interest due to Zariski-Nagata Theorem \cite[Theorem 3.14]{Eisenbud}.

\begin{theorem} (Zariski-Nagata)
	Let $I$ be a radical homogeneous  ideal, and let char$( \K)=0$. For $m\geq 1$
	$$I^{(m)}=\{f \in \K[x_0,\ldots, x_N]\; : \; \text{f vanishes to order }\geq \text{m in all points } P \in V(I)\} .$$
\end{theorem}
In the situation when  $V(I)$ is a finite set of points $P_1,\ldots,P_t \in \P^N$, the symbolic power is particularly easy to compute: 
$$I^{(m)} = \bigcap_{i=1}^{t} I(P_i)^m. $$

An arrangements of lines $\mathcal{A}$ is a finite set of mutually distinct lines 
$L_1, \ldots, L_s$. An arrangement of lines determines a finite set of points 
$\mathcal{P}=\{P_1,\ldots, P_t\}$ in $\P^2$, where at least $2$ of arrangement lines intersect.
For $i\geq 2$, we denote by $t_i(\mathcal{A})$ the number of points in $\mathcal{P}$ where exactly $i$ lines from $\mathcal{A}$ intersect. These numbers define the $t$-vector $t(\mathcal{A})$ of $\mathcal{A}$
\[t(\mathcal{A})=(t_2(\mathcal{A}),\:t_3(\mathcal{A}),\ldots,\:t_s(\mathcal{A})).\]
It is a basic combinatorial invariant of $\mathcal{A}$.

For line arrangements defined over $\R$, the following property has been distinguished. 
 
\begin{definition}[Simplicial arrangement]
 We say that an arrangement $\mathcal{A}$ of real lines is simplicial  if every connected component of its complement $\P^2(\R)\setminus \mathcal{A}$ is a triangle. 
\end{definition}

It is expected, but not known, if (apart of $3$ obvious infinite families described in \cite{Gru}) there are only finitely many sporadic examples. A list of such examples was constructed by Gr\"unbaum in \cite{Gru} and extended recently by Cuntz in \cite{Cuntz}.

\section{Simplicial arrangements $\mathcal{A}(31,2)$ and $\mathcal{A}(31,3)$}
\label{sec:A312andA313}

The arrangements we study here come from \cite{Gru}, where they are called $\mathcal{A}(31,2)$ and $\mathcal{A}(31,3)$.

Configurations $\mathcal{A}(31,2)$ and $\mathcal{A}(31,3)$ are non isomorphic simplicial arrangements of $31$ lines with the total number of $127$ intersection points. Moreover we have
$$t_2 = 54, t_3 =42, t_4= 21, t_5=6, t_6=1, t_8=3$$
and all other $t_i =0.$

\subsection{Configuration $\mathcal{A}(31,3)$} 
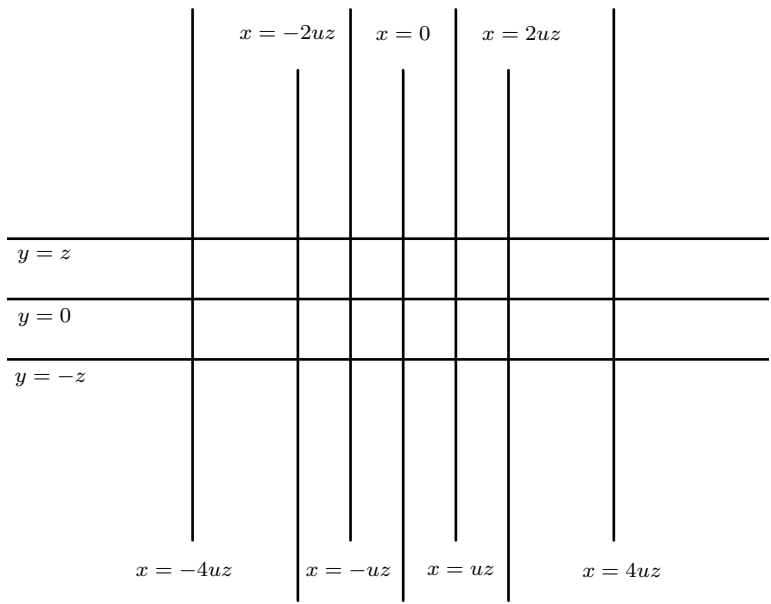
\begin{figure}[H] 
	\centering
	\begin{tikzpicture}
	[line cap=round,line join=round,>=triangle 45,x=1.0cm,y=1.0cm, scale=0.8]
	\clip(-6.5,-6) rectangle (6,5);
	\draw [line width=1.pt] (0.8660254037844386,-4) -- (0.8660254037844386,4.8);
	\draw [line width=1.pt] (1.7320508075688772,-5) -- (1.7320508075688772,3.792672179572946);
	\draw [line width=1.pt] (3.4641016151377544,-4) -- (3.4641016151377544,4.8);
	\draw [line width=1.pt] (-0.8660254037844386,-4) -- (-0.8660254037844386,4.8);
	\draw [line width=1.pt] (0.,-5) -- (0.,3.792672179572946);
	\draw [line width=1.pt] (-1.7320508075688772,-5) -- (-1.7320508075688772,3.792672179572946);
	\draw [line width=1.pt] (-3.4641016151377544,-4) -- (-3.4641016151377544,4.8);
	\draw [line width=1.pt,domain=-10:11.621005370643537] plot(\x,{(--1.-0.*\x)/1.});
	\draw [line width=1.pt,domain=-10:11.621005370643537] plot(\x,{(-0.-0.*\x)/1.});
	\draw [line width=1.pt,domain=-10:11.621005370643537] plot(\x,{(-1.-0.*\x)/1.});
	\begin{scriptsize}
	\draw[color=black] (0,4.4) node {$x = 0$};
	\draw[color=black] (-0.9,-4.5) node {$x = -uz$};
	\draw[color=black] (0.95,-4.5) node {$x = uz$};
	\draw[color=black] (-1.9,4.4) node {$x = -2uz$};
	\draw[color=black] (1.95,4.4) node {$x = 2uz$};
	\draw[color=black] (-3.6,-4.5) node {$x = -4uz$};
	\draw[color=black] (3.6,-4.5) node {$x = 4uz$};
	\draw[color=black] (-5.9,-0.3) node {$y = 0$};
	\draw[color=black] (-5.9,0.7) node {$y = z$};
	\draw[color=black] (-5.8,-1.3) node {$y = -z$};
	\end{scriptsize}
	\end{tikzpicture}
	\caption{. $10$ initial lines in the construction of $\mathcal{A}(31,3)$.}
	\label{fig:tenlines}
\end{figure}

This configuration can be realized in the following way. We begin with ten lines:
$$x \pm auz=0, \; y\pm bz=0,$$ where $u= \frac{\sqrt{3}}{2}$, $a=0,1,2,4$ and $b=0,1$. These ten lines are visualized in Figure \ref{fig:tenlines}. 

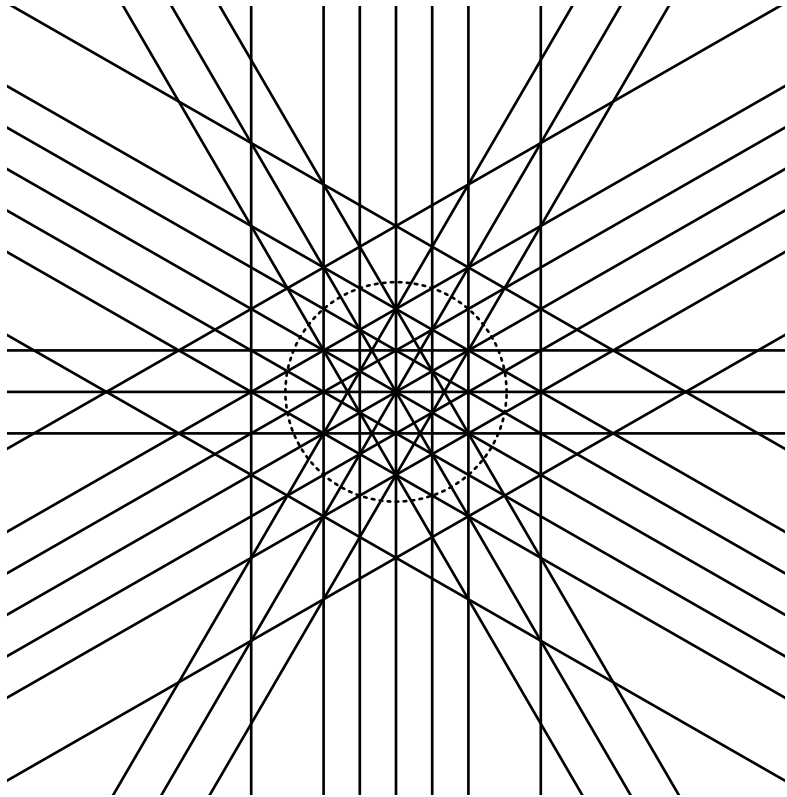
\begin{figure}[h]
	\centering
	\begin{tikzpicture}[line cap=round,line join=round,>=triangle 45,x=1.0cm,y=1.0cm,scale=0.55]
	\clip(-9.3,-9.7) rectangle (9.3,9.3);
	\draw [line width=1.pt] (0.8660254037844386,-10.003080137831516) -- (0.8660254037844386,9.58127622309531);
	\draw [line width=1.pt] (1.7320508075688772,-10.003080137831516) -- (1.7320508075688772,9.58127622309531);
	\draw [line width=1.pt] (3.4641016151377544,-10.003080137831516) -- (3.4641016151377544,9.58127622309531);
	\draw [line width=1.pt] (-0.8660254037844386,-10.003080137831516) -- (-0.8660254037844386,9.58127622309531);
	\draw [line width=1.pt] (0.,-10.003080137831516) -- (0.,9.58127622309531);
	\draw [line width=1.pt] (-1.7320508075688772,-10.003080137831516) -- (-1.7320508075688772,9.58127622309531);
	\draw [line width=1.pt] (-3.4641016151377544,-10.003080137831516) -- (-3.4641016151377544,9.58127622309531);
	\draw [line width=1.pt,domain=-10.669652404943239:10.026086142109872] plot(\x,{(--1.-0.*\x)/1.});
	\draw [line width=1.pt,domain=-10.669652404943239:10.026086142109872] plot(\x,{(-0.-0.*\x)/1.});
	\draw [line width=1.pt,domain=-10.669652404943239:10.026086142109872] plot(\x,{(-1.-0.*\x)/1.});
	\draw [line width=1.pt,domain=-10.669652404943239:10.026086142109872] plot(\x,{(-0.--0.8660254037844386*\x)/0.5});
	\draw [line width=1.pt,domain=-10.669652404943239:10.026086142109872] plot(\x,{(--1.--0.8660254037844386*\x)/0.5});
	\draw [line width=1.pt,domain=-10.669652404943239:10.026086142109872] plot(\x,{(-1.--0.8660254037844386*\x)/0.5});
	\draw [line width=1.pt,domain=-10.669652404943239:10.026086142109872] plot(\x,{(--3.4641016151377544-0.5*\x)/0.8660254037844386});
	\draw [line width=1.pt,domain=-10.669652404943239:10.026086142109872] plot(\x,{(--1.7320508075688772-0.5*\x)/0.8660254037844386});
	\draw [line width=1.pt,domain=-10.669652404943239:10.026086142109872] plot(\x,{(--0.8660254037844386-0.5*\x)/0.8660254037844386});
	\draw [line width=1.pt,domain=-10.669652404943239:10.026086142109872] plot(\x,{(-0.-0.5*\x)/0.8660254037844386});
	\draw [line width=1.pt,domain=-10.669652404943239:10.026086142109872] plot(\x,{(-0.8660254037844386-0.5*\x)/0.8660254037844386});
	\draw [line width=1.pt,domain=-10.669652404943239:10.026086142109872] plot(\x,{(-1.7320508075688772-0.5*\x)/0.8660254037844386});
	\draw [line width=1.pt,domain=-10.669652404943239:10.026086142109872] plot(\x,{(-3.4641016151377544-0.5*\x)/0.8660254037844386});
	\draw [line width=1.pt,domain=-10.669652404943239:10.026086142109872] plot(\x,{(-1.--0.8660254037844388*\x)/-0.5});
	\draw [line width=1.pt,domain=-10.669652404943239:10.026086142109872] plot(\x,{(-0.--0.8660254037844388*\x)/-0.5});
	\draw [line width=1.pt,domain=-10.669652404943239:10.026086142109872] plot(\x,{(--1.--0.8660254037844388*\x)/-0.5});
	\draw [line width=1.pt,domain=-10.669652404943239:10.026086142109872] plot(\x,{(--3.4641016151377544--0.5*\x)/0.8660254037844388});
	\draw [line width=1.pt,domain=-10.669652404943239:10.026086142109872] plot(\x,{(--1.7320508075688772--0.5*\x)/0.8660254037844388});
	\draw [line width=1.pt,domain=-10.669652404943239:10.026086142109872] plot(\x,{(--0.8660254037844386--0.5*\x)/0.8660254037844388});
	\draw [line width=1.pt,domain=-10.669652404943239:10.026086142109872] plot(\x,{(-0.--0.5*\x)/0.8660254037844388});
	\draw [line width=1.pt,domain=-10.669652404943239:10.026086142109872] plot(\x,{(-0.8660254037844386--0.5*\x)/0.8660254037844388});
	\draw [line width=1.pt,domain=-10.669652404943239:10.026086142109872] plot(\x,{(-1.7320508075688772--0.5*\x)/0.8660254037844388});
	\draw [line width=1.pt,domain=-10.669652404943239:10.026086142109872] plot(\x,{(-3.4641016151377544--0.5*\x)/0.8660254037844388});
	\draw [line width=1.pt,dotted] (0,0) circle (2.64575cm);
	\end{tikzpicture}
	\caption{. Configuration $A(31,3)$. The line $z=0$, which lies at infinity, is not shown. Intersection points of lines which belong to the dotted circle form the orbit represented by point $(9:2u:4u)$.}
	\label{fig:A31_3}
\end{figure}

Then we rotate these lines by $60^{\circ}$ and $120^{\circ}$ around the point $(0:0:1)$. In this way, we obtain $30$ lines. The last line is the line at infinity $z=0$. As a result we obtain a configuration of lines  indicated in Figure \ref{fig:A31_3}. Taking the product of linear forms defining the ten initial lines we obtain the following polynomial
$$F_{10}=x^7y^3-x^7yz^2-\frac{63}{4}x^5y^3z^2+\frac{63}{4}x^5yz^4+
\frac{189}{4}x^3y^3z^4-\frac{189}{4}x^3yz^6-27xy^3z^6+27xyz^8$$
and taking the product of all $31$ lines we get a polynomial $F_{31}$ of degree $31$. We are interested in the Jacobian ideal $Jac(F_{31})$ defined by this polynomial. The radical $I_{\mathcal{A}(31,3)}$ of this ideal describes all $127$ intersection points among arrangement lines. 

By construction the arrangement  is invariant under the group $G = \mathbb{Z}_{3} \times \mathbb{Z}_{2}$ generated by the rotation matrix $ A= \begin{bmatrix}
\frac{1}{2} & -u & 0  \\
u & \frac{1}{2} & 0 \\
0 & 0 & 1
\end{bmatrix}$
and the reflection $ B= \begin{bmatrix}
-1 & 0 & 0  \\
0 & 1 & 0 \\
0 & 0 & 1
\end{bmatrix}$, which together generate the dihedral group $D_{6}$. This group acts on the set of $127$ intersection points so that, there are the following orbits

{
	\begin{center}
		\begin{tabular}[H!]{c|c|c}
			length of orbit & number of orbits & representing point\\
			\hline
			\rule{0pt}{3ex} 
			$1$	& $1$ & $(0:0:1)$\\
			\rule{0pt}{4ex} 
			
			$6$	& $9$ &  $(0:1:1)$, $(1:0:u)$, $(2u:0:1)$, \\
			&  & $(0:2:1)$, $(4u:0:1)$, $(0:4:1)$, \\
			& & $(u:0:1)$, $(8u:0:1)$, $(1:0:1)$\\
			
			\rule{0pt}{4ex}
			$12$	& $6$ &  $(6u:1:4)$, $(9:2u:4u)$, $(4u:1:1)$,\\
			& & $(15:6u:4u)$, $(6u:1:1)$, $(10u:1:1)$

		\end{tabular}
		
	\end{center}
}
It is helpful to consider the sub-arrangement $\mathcal{B}_{21}$ consisting of $7$ lines:
$$x\pm auz=0, \; y \pm bz=0,$$
where $a=1,4$, $b=0,1$ and images of these lines under $A$ and $A^2$. This $21$ lines intersect altogether in $115$ points, with multiplicities $t_2 =72$, $t_3 =40$ and $t_4 =3$. The difference between the $127$ and $115$ points is one full orbit represented by point  $(9:2u:4u)$.

The points in this orbit are now contained each in only one of the $21$ lines. In order to get an element in $I^{(3)}_{\mathcal{A}(31,3)}$ we need to complete the $21$ lines by a divisor $\Gamma$ vanishing in these $12$ points to order $2$ and passing through the remaining $72$ points, which are double points for $\mathcal{B}_{21}$.

To this end we consider $X=\mathbb{P}^{2} / G$. The ring of invariant polynomials $\mathbb{K}[x,y,z]^{G}$ is generated by 
$$f_1= z, \quad f_2 = x^2 + y^2, \quad f_3 = 11x^6+15x^4 y^2 +45x^2y^4 +9y^6.$$
Using Moliens's Theorem (see \cite{Sturm}, Theorem $2.2.1$), we see that the space of invariant polynomials of degree $12$ has dimension $12$.

Since vanishing to order $2$ at a smooth point of $X$ imposes $3$ conditions and the $72$ points split into $4$ orbits of order $6$ and $4$ orbits of order $12$, counting conditions 
$$12-3-4-4=1>0$$
we conclude that the desired divisor $\Gamma$ exists (it is invariant under $G$, so it pulls back from $X$). Computing with Singular, we are able to identify the equation of $\Gamma$:

\begin{align*}
\frac{2093688}{17}f_1^{12}-\frac{9398511}{34}f_1^{10} f_2 +\frac{2995218}{17} f_1^8 f_2^2-\frac{64485153}{1088} f_1^6  f_2^3 \\ + \frac{18708003}{4352}f_1^4  f_2 ^4  + \frac{1258659}{4352}f_1^2  f_2^5- \frac{493695}{4352}f_2^6+ \frac{2121309}{1088}f_1^6  f_3 \\ - \frac{402561}{4352} f_1^4  f_2 f_3 - \frac{158697}{4352}f_1^2  f_2^2  f_3  + \frac{2619}{128} f_2^3   f_3- \frac{3979}{4352} f_3^2,
\end{align*}
in terms of the invariant generators $f_1$, $f_2$, $f_3$.

Considering the equation of $\Gamma$ in the ring $\mathbb{K}[f_1, f_2, f_3]$, it is easy to check that there is just one singular point, which is locally simple crossing. This implies that $\Gamma$ is irreducible.

For the non-containment  $I^{(3)}_{\mathcal{A}(31,3)} \nsubseteq I^2_{\mathcal{A}(31,3)}$ we used Singular. We do not have a theoretical proof. Summing up claims in this section, we see that Theorem A is proved.

\begin{figure}
	\centering
	\includegraphics[height=10cm,width=14cm]{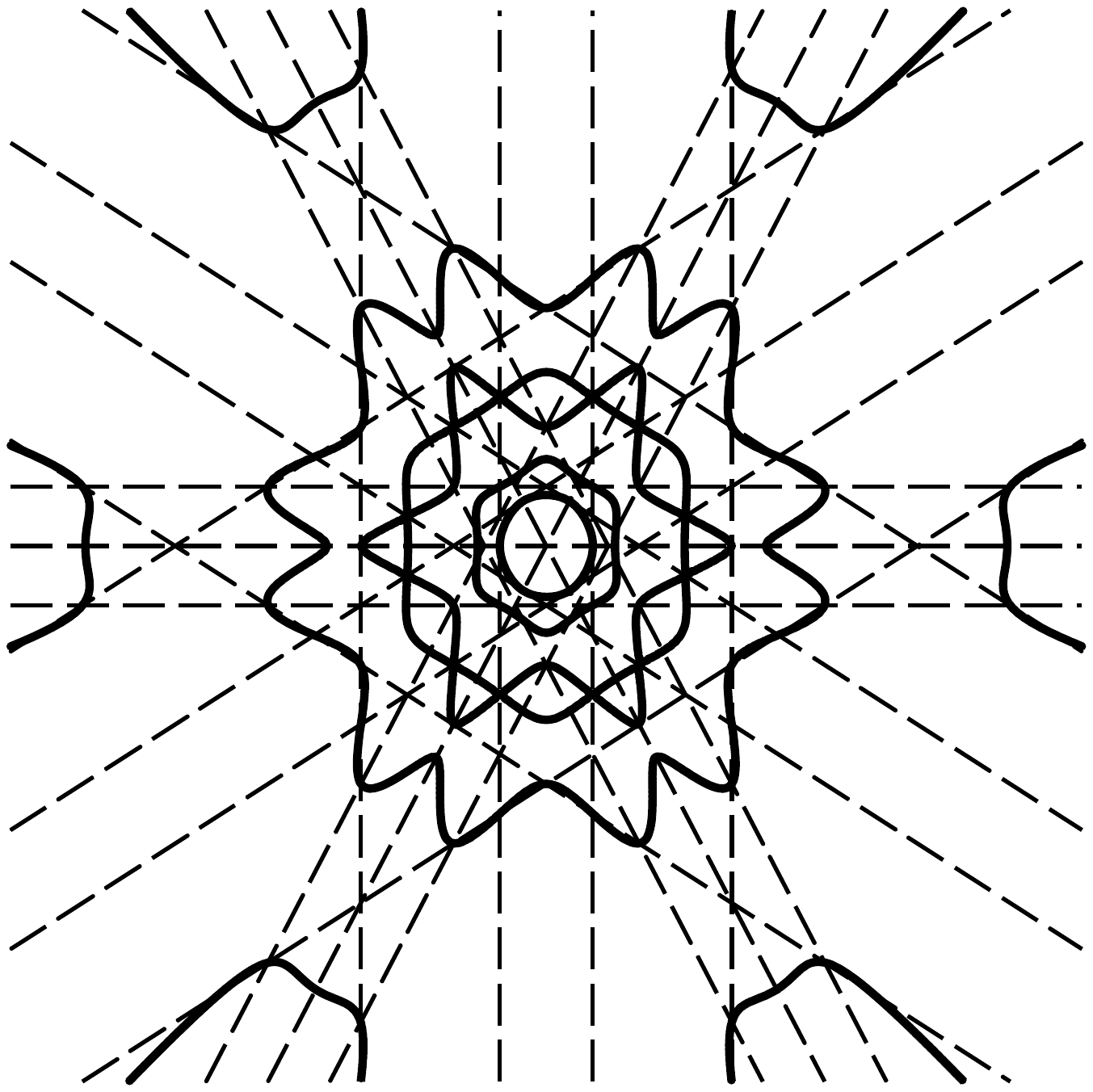}
	\caption{. The graph of affine part of polynomial from set $I_{\mathcal{A}(31,3)} ^{(3)}\setminus I_{\mathcal{A}(31,3)}^2$, which consists of $21$ dashed lines and a curve of degree $12$.}
	\label{fig:high deg curve and lines}
\end{figure}

At the end of this section we want to underline another interesting observation about curve from set $I_{\mathcal{A}(31,3)} ^{(3)}\setminus I_{\mathcal{A}(31,3)}^2$ indicated on Figure \ref{fig:high deg curve and lines}. The twelve visible double points are the only singular points for this curve, thus we can easily calculate the arithmetic genus, which is 

$$g=\binom{12-1}{2}-12\binom{2}{2}= 43.$$ 

This is the first known example of the curve, which form an element from the set $I^{(3)} \setminus I^2$ and which is not rational at the same time.

\subsection{Configuration $\mathcal{A}(31,2)$}

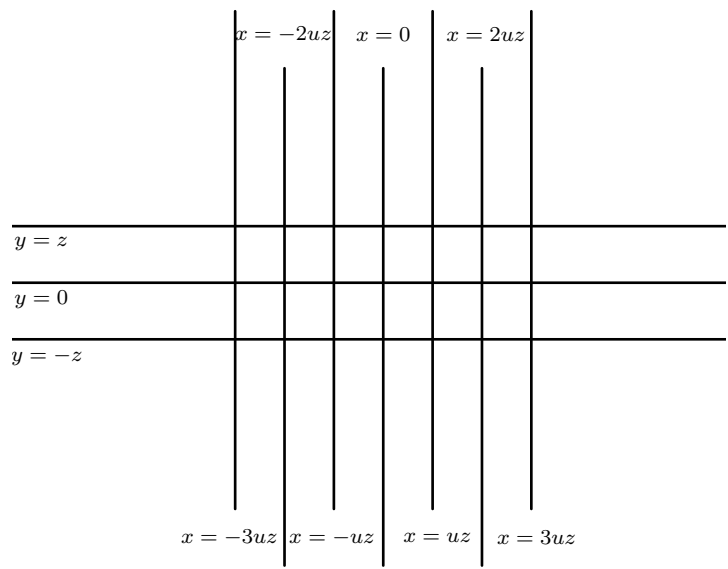
\begin{figure}
	\centering
	\begin{tikzpicture}
	[line cap=round,line join=round,>=triangle 45,x=1.0cm,y=1.0cm, scale=0.75]
	\clip(-6.5,-6) rectangle (6,5);
	\draw [line width=1.pt] (0.8660254037844386,-4) -- (0.8660254037844386,4.8);
	\draw [line width=1.pt] (1.7320508075688772,-5) -- (1.7320508075688772,3.792672179572946);
	\draw [line width=1.pt] (2.598,-4) -- (2.598,4.8);
	\draw [line width=1.pt] (-0.8660254037844386,-4) -- (-0.8660254037844386,4.8);
	\draw [line width=1.pt] (0.,-5) -- (0.,3.792672179572946);
	\draw [line width=1.pt] (-1.7320508075688772,-5) -- (-1.7320508075688772,3.792672179572946);
	\draw [line width=1.pt] (-2.598,-4) -- (-2.598,4.8);
	\draw [line width=1.pt,domain=-10:11.621005370643537] plot(\x,{(--1.-0.*\x)/1.});
	\draw [line width=1.pt,domain=-10:11.621005370643537] plot(\x,{(-0.-0.*\x)/1.});
	\draw [line width=1.pt,domain=-10:11.621005370643537] plot(\x,{(-1.-0.*\x)/1.});
	\begin{scriptsize}
	\draw[color=black] (0,4.4) node {$x = 0$};
	\draw[color=black] (-0.9,-4.5) node {$x = -uz$};
	\draw[color=black] (0.95,-4.5) node {$x = uz$};
	\draw[color=black] (-1.75,4.4) node {$x = -2uz$};
	\draw[color=black] (1.8,4.4) node {$x = 2uz$};
	\draw[color=black] (-2.7,-4.5) node {$x = -3uz$};
	\draw[color=black] (2.69,-4.5) node {$x = 3uz$};
	\draw[color=black] (-6,-0.3) node {$y = 0$};
	\draw[color=black] (-6,0.7) node {$y = z$};
	\draw[color=black] (-5.9,-1.3) node {$y = -z$};
	\end{scriptsize}
	\end{tikzpicture}
	\caption{. $10$ initial lines in the construction of $\mathcal{A}(31,2)$.}
	\label{fig:A312_tenLines}
\end{figure}

This configuration is very similar to $\mathcal{A}(31,3)$. It can be realized starting with lines
$$x \pm auz=0, \; y\pm bz=0,$$ where  $a=0,1,2,3$ and $b=0,1$. These lines are visualized in Figure \ref{fig:A312_tenLines}.

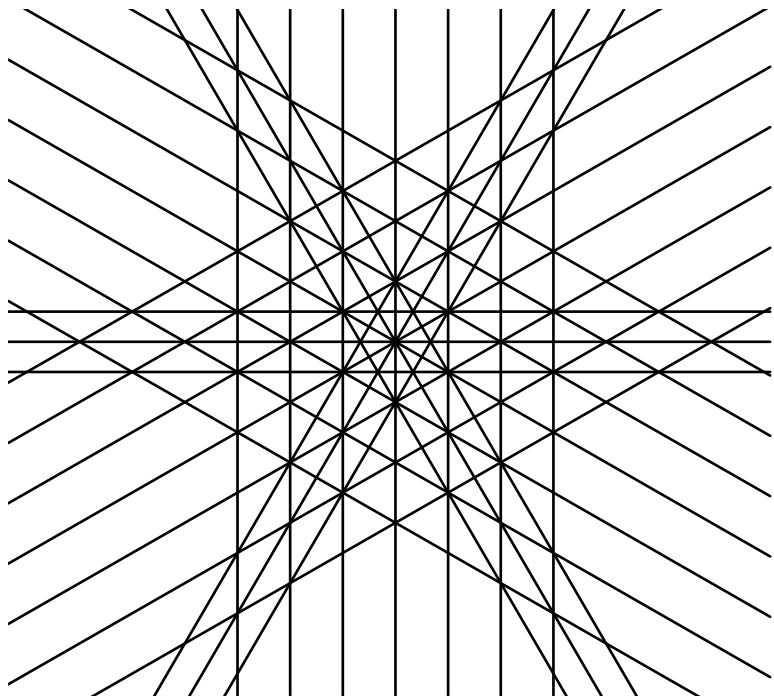
\begin{figure}
	\centering
	\begin{tikzpicture}[line cap=round,line join=round,>=triangle 45,x=1.0cm,y=1.0cm,scale=0.4]
	\clip(-12.73374002469864,-11.717725054379748) rectangle (12.9,11.01595705783588);
	\draw [line width=1.pt] (1.7320508075688772,-11.717725054379748) -- (1.7320508075688772,11.01595705783588);
	\draw [line width=1.pt] (3.4641016151377544,-11.717725054379748) -- (3.4641016151377544,11.01595705783588);
	\draw [line width=1.pt] (5.196152422706632,-11.717725054379748) -- (5.196152422706632,11.01595705783588);
	\draw [line width=1.pt] (-1.7320508075688772,-11.717725054379748) -- (-1.7320508075688772,11.01595705783588);
	\draw [line width=1.pt] (-3.4641016151377544,-11.717725054379748) -- (-3.4641016151377544,11.01595705783588);
	\draw [line width=1.pt] (0.,-11.717725054379748) -- (0.,11.01595705783588);
	\draw [line width=1.pt] (-5.196152422706632,-11.717725054379748) -- (-5.196152422706632,11.01595705783588);
	\draw [line width=1.pt,domain=-12.73374002469864:12.326711565445072] plot(\x,{(-0.-0.*\x)/1.});
	\draw [line width=1.pt,domain=-12.73374002469864:12.326711565445072] plot(\x,{(-1.-0.*\x)/1.});
	\draw [line width=1.pt,domain=-12.73374002469864:12.326711565445072] plot(\x,{(--1.-0.*\x)/1.});
	\draw [line width=1.pt,domain=-12.73374002469864:12.326711565445072] plot(\x,{(--1.7320508075688772-0.5*\x)/0.8660254037844386});
	\draw [line width=1.pt,domain=-12.73374002469864:12.326711565445072] plot(\x,{(-1.--0.8660254037844386*\x)/0.5});
	\draw [line width=1.pt,domain=-12.73374002469864:12.326711565445072] plot(\x,{(-0.--0.8660254037844386*\x)/0.5});
	\draw [line width=1.pt,domain=-12.73374002469864:12.326711565445072] plot(\x,{(--1.--0.8660254037844386*\x)/0.5});
	\draw [line width=1.pt,domain=-12.73374002469864:12.326711565445072] plot(\x,{(--3.4641016151377544-0.5*\x)/0.8660254037844386});
	\draw [line width=1.pt,domain=-12.73374002469864:12.326711565445072] plot(\x,{(--5.196152422706632-0.5*\x)/0.8660254037844386});
	\draw [line width=1.pt,domain=-12.73374002469864:12.326711565445072] plot(\x,{(-1.7320508075688772-0.5*\x)/0.8660254037844386});
	\draw [line width=1.pt,domain=-12.73374002469864:12.326711565445072] plot(\x,{(-3.4641016151377544-0.5*\x)/0.8660254037844386});
	\draw [line width=1.pt,domain=-12.73374002469864:12.326711565445072] plot(\x,{(-5.196152422706632-0.5*\x)/0.8660254037844386});
	\draw [line width=1.pt,domain=-12.73374002469864:12.326711565445072] plot(\x,{(-0.-0.5*\x)/0.8660254037844386});
	\draw [line width=1.pt,domain=-12.73374002469864:12.326711565445072] plot(\x,{(-0.--0.5*\x)/0.8660254037844388});
	\draw [line width=1.pt,domain=-12.73374002469864:12.326711565445072] plot(\x,{(--1.7320508075688772--0.5*\x)/0.8660254037844388});
	\draw [line width=1.pt,domain=-12.73374002469864:12.326711565445072] plot(\x,{(--3.4641016151377544--0.5*\x)/0.8660254037844388});
	\draw [line width=1.pt,domain=-12.73374002469864:12.326711565445072] plot(\x,{(--5.196152422706632--0.5*\x)/0.8660254037844388});
	\draw [line width=1.pt,domain=-12.73374002469864:12.326711565445072] plot(\x,{(-1.7320508075688772--0.5*\x)/0.8660254037844388});
	\draw [line width=1.pt,domain=-12.73374002469864:12.326711565445072] plot(\x,{(-3.4641016151377544--0.5*\x)/0.8660254037844388});
	\draw [line width=1.pt,domain=-12.73374002469864:12.326711565445072] plot(\x,{(-5.196152422706632--0.5*\x)/0.8660254037844388});
	\draw [line width=1.pt,domain=-12.73374002469864:12.326711565445072] plot(\x,{(--1.--0.8660254037844388*\x)/-0.5});
	\draw [line width=1.pt,domain=-12.73374002469864:12.326711565445072] plot(\x,{(-0.--0.8660254037844388*\x)/-0.5});
	\draw [line width=1.pt,domain=-12.73374002469864:12.326711565445072] plot(\x,{(-1.--0.8660254037844388*\x)/-0.5});
	\end{tikzpicture}
	\caption{. Configuration $\mathcal{A}(31,2)$. The line $z=0$ at infinity is not shown.}
	\label{fig:A312}
\end{figure}

Rotating again by $60^{\circ}$ and $120^{\circ}$ and taking the line at infinity we obtain the configuration presented in Figure \ref{fig:A312}.

The multiplicities vector of this configuration is the same as for $\mathcal{A}(31,3)$, i.e., there is 
$$t_2 = 54, t_3 =42, t_4= 21, t_5=6, t_6=1, t_8=3.$$
In particular there are again $127$ intersection points of pairs of arrangement lines. However, a quick Singular check shows that now we have 
$$I^{(3)}_{\mathcal{A}(31,2)} \subseteq I^{2}_{\mathcal{A}(31,2)} .$$
This shows, once again (see \cite{combinatorics}), that the (non)containment property is quite subtle and cannot be decided by looking at the basis combinatorial invariants only.

\subsection{Arrangement $\mathcal{B}_{21}$}

Now we consider more closely the arrangement $\mathcal{B}_{21}$ defined in the previous section. We keep the notation introduced there.

The ideal $I_{\mathcal{B}_{21}}$ defines $115$ points. This is the subset of $127$ points defined by $I_{\mathcal{A}(31,3)}$, the difference being one $G-$orbit, consisting of $12$ double points of $\Gamma$. 

In order to exhibit an element in $I^{(3)}_{\mathcal{B}_{21}}$, we need to find a divisor vanishing at the $72$ points, where only $2$ of arrangement lines meet. 

Revoking again Molien's Theorem, we see that the dimension of the space of $G-$invariant polynomials of degree $10$ is $9$, thus the expected dimension of invariant polynomials vanishing at $4$ order $6$ and $4$ order $12$ orbits in which the $72$ points split is 
$$9-4-4=1>0.$$
Hence there is a divisor $\Delta$ of degree $10$ vanishing at these $72$ points. We can express its equation in terms of invariant polynomials:
\begin{align*}
-\frac{38320128}{107}f_1^{10}+ \frac{80453952}{107}f_1^8 f_2 -\frac{42393996}{4107}f_1^6 f_2^2 + \frac{50759217}{214} f_1^4 f_2^3 \\  -\frac{20519091}{856} f_1^2 f_2^4+ \frac{67086}{107} f_2^5- \frac{3811059}{214} f_1^4 f_3 + \frac{1778227}{856} f_1^2 f_2 f_3 - \frac{6089}{107} f_2^2 f_3.  \end{align*}
Since $\Delta$ is smooth, it is irreducible.

The  non-containment $I^{(3)}_{\mathcal{B}_{21}} \nsubseteq I^{2}_{\mathcal{B}_{21}}$ is proved again with the aid of Singular \cite{DGPS}. 

Summing up the claims of this section, we obtain the proof of Theorem B.

\begin{figure}[H]
	\centering
	\includegraphics[height=10cm,width=14cm]{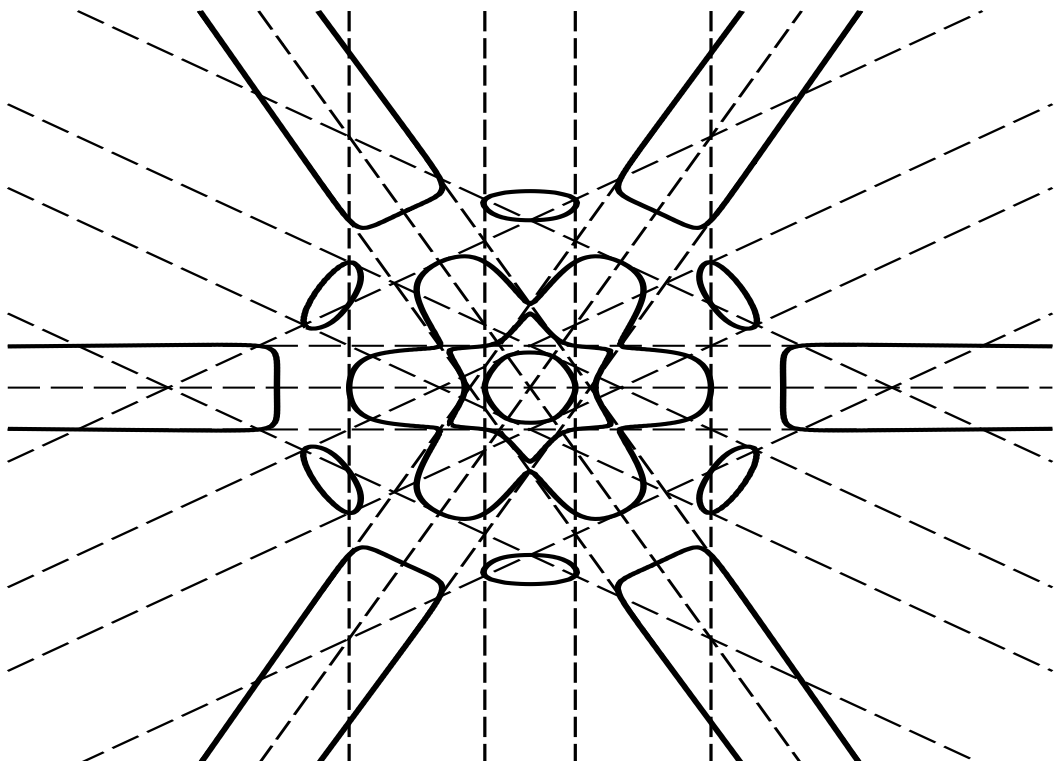}
	\caption{. The graph of affine part of polynomial from set $I_{\mathcal{B}_{21}}^{(3)}\setminus I_{\mathcal{B}_{21}}^2$, which consists of $21$ dashed lines and a curve of degree $10$.}
	\label{fig:F10}
\end{figure}

As for the curve in Section \ref{sec:A312andA313}, we also calculate arithmetic genus for the curve, which is indicated as a solid line on Figure \ref{fig:F10}. Using any computer algebra program one can check that this curve has only four non-reduced singular points, and that its genus is $g=30.$ This means that this curve is not rational.

We conclude this section by noting that the arrangement $\mathcal{B}_{21}$ is not free. Indeed, its characteristic polynomial is 
$$\chi(\mathcal{B}_{21},t)=1-20t+141t^2 $$ and it does not split over the integers (see Main Theorem in  \cite{terao}).

\section{Realizability over rational numbers}
The first non-containment 
$$I^{(3)} \nsubseteq I^2$$
was the dual Hesse arrangement, see \cite{DST13}. This arrangement cannot be realized over the reals. The first real non-containment example, the B\"or\"oczky arrangement of $12$ lines was discovered in \cite{CGMLLPS2015}. It was realized in \cite{BoHa} and \cite{MJ2015} that the B\"or\"oczky arrangement can be defined over the rational numbers. Additional examples were provided in \cite{malara-szpond} and \cite{combinatorics}. Such examples are quite rare, so we find it worth to mention that $\mathcal{A}(31,3)$ and $\mathcal{B}_{21}$ can be both realized over $\mathbb{Q}$. We can be quite explicit here. Table \ref{lines} contains equations of all $31$ lines, whereas coordinates of their intersection points are provided in Table \ref{points}.

{ \renewcommand{\arraystretch}{1.5}\begin{table}[h]
		\begin{center}
			\begin{tabular}{l l} 
				\multicolumn{2}{c}{$\mathcal{A}(31,3)$} \\
				\hline
				
				$x+y+iz=0$, & for $i \in \{0,2,3,4,5,6,8\}$,  \\
				
				$2x-y+jz=0$, & for $j \in \{4,6,7,8,9,10,12\}$, \\
				
				$3x+kz=0$, &  for $k \in \{8,10,11,12,13,14,16\}$,\\
				
				$x-2y+lz=0$, & for $l \in \{2,4,6\}$, \\

 				$4x+y+mz=0$, & for $m \in \{14,16,18\}$,\\
				
				$5x-y+nz=0$, & for $n \in \{18,20,22\}$,\\
				
				$z=0$
			\end{tabular}
			
		\end{center}
		\caption{$\:$ The equations of lines of $\mathcal{A}(31,3)$.}
		\label{lines}
	\end{table}
}

{
	\renewcommand{\arraystretch}{1.3} 
	\renewcommand{\tabcolsep}{2.5pt}
	 \begin{table}[H]
	\begin{center}
		\begin{tabular}{c|cccccc}
			\multicolumn{7}{c}{$\mathcal{A}(31,3)$} \\
			\hline
			\multirow{8}{*}{double} 
			 
			&  (2,-2,-3),& (6,-6,-1), &(7,-7,-3),& (13,-13,-3),&
			   (3,-1,-1),& (13,-7,-3),\\

			&   (5,1,-1),& (11,7,-3),& (22,2,-3),& (2,6,-1),&		
			    (17,7,-3),& (11,13,-3),\\

			&    (7,2,-3),& (3,0,-1),& (14,-5,-3),& (16,-7,-3),&
			     (7,-1,-2),& (23,-5,-6),\\

			&    (25,-7,-6),& (5,0,-1),& (17,-2,-3),& (8,7,-3),& 
			     (10,5,-3),& (23,7,-6),\\

			&   (25,5,-6),& (9,1,-2),& (-22,-8,3),& (-22,-5,6),& 
			    (-26,-7,6),& (-22,1,6),\\

			&   (-26,-1,6),& (-2,8,3),& (-22,7,6),& (-26,5,6),&
			    (-6,-8,1),& (-13,-14,3),\\

			& (-13,-8,3),& (-11,8,3),& (-2,8,1),& (-11,14,3),&
			  (-8,-22,3),& (8,-11,-3),\\	
			
			& (8,-26,-3),& (10,-7,-3),& (14,7,-3),& (-16,22,3),&
			  (16,11,-3),& (16,26,-3), \\
			
			&  (7,0,-2),& (23,-8,-6),& (23,4,-6),& (25,-4,-6),&
			   (25,8,-6),& (9,0,-2)\\
			\hline
			\multirow{7}{*}{triple} 
				
			& (4,-4,-3),& (2,-2,-1),& (8,-8,-3),& (4,-4,-1),& 
			(14,-14,-3),& (16,-16,-3),\\

			&  (3,-3,-1),& (11,-11,-3),& (2,0,-1),& (16,-10,-3),&
			   (8,4,-3),&  (16,-4,-3),\\

			&   (6,0,-1),& (8,10,-3),& (6,2,-1),& (20,4,-3),& 
			    (4,4,-1),& (16,8,-3),\\

			&    (8,16,-3),& (10,14,-3),& (5,3,-1),& (13,11,-3),&
			     (8,1,-3),&  (5,-2,-1),\\

			&  (16,-1,-3),& (3,2,-1),& (-16,-5,3),& (-34,-8,9),& 
			   (-38,-10,9),& (-32,2,9),\\

			&  (-40,-2,9),& (-8,5,3),& (-34,10,9),& (-38,8,9),&		
			   (-14,-16,3),& (-16,-20,3),\\

			&  (-11,-10,3),& (-16,-14,3),& (-8,14,3),& (-8,20,3),&
			   (-10,16,3),& (-13,10,3)\\	
			\hline
			\multirow{4}{*}{quadruple}
			    
			& (2,1,0),& (-1,4,0),& (1,5,0),& (10,-10,-3),& 
			  (8,-2,-3),& (14,-8,-3),\\

			&  (11,-5,-3),& (11,1,-3),& (13,-1,-3),& (16,2,-3),& 
			(10,8,-3),& (13,5,-3),\\

			& (14,10,-3),& (10,-1,-3),& (4,-1,-1),& (11,-2,-3),&
			  (13,-4,-3),& (4,1,-1),\\

			& (14,1,-3),& (11,4,-3),& (13,2,-3)& & & \\
			\hline
			\multirow{1}{*}{quintuple} 
			    
			&  (10,-4,-3),& (4,-2,-1),& (10,2,-3),& (14,-2,-3),& 		               (14,4,-3),& (4,2,-1)\\
			\hline
			\multirow{1}{*}{sextuple} 
			 
			&  (4,0,-1)\\
			\hline
			\multirow{1}{*}{octuple} 
			 
			&   (-1,1,0),& (1,2,0),& (0,1,0)& & & \\
			
		\end{tabular}
		\caption{$\:$ The coordinates of points of $\mathcal{A}(31,3)$.}
	\label{points}
		
	\end{center} \end{table}
}
\paragraph*{Acknowledgements.}
We would like to warmly thank T. Szemberg for all helpful remarks, valuable comments and inspiring discussions  which greatly improve the original draft of our paper.

The research of Lampa-Baczy\'nska was partially supported by National Science Centre, Poland, grant 2016/23/N/ST1/01363, the research of Malara was partially supported by National Science Centre, Poland, grant 2016/21/N/ST1/01491.


\bigskip \footnotesize

\bigskip
\noindent Marek Janasz, Magdalena~Lampa-Baczy\'nska, Grzegorz Malara \\
   Instytut Matematyki UP,
   Podchor\c a\.zych 2,
   PL-30-084 Krak\'ow, Poland
\\
\nopagebreak
  \textit{E-mail address:} \texttt{mjanasz@op.pl}\\
  \textit{E-mail address:} \texttt{lampa.baczynska@wp.pl}\\
  \textit{E-mail address:} \texttt{grzegorzmalara@gmail.com, gmalara@impan.pl}
  
\vspace{0.5cm}
\noindent
Grzegorz Malara current address: Institute of Mathematics, Polish Academy of Sciences, \\
\'Sniadeckich 8, PL-00-656 Warszawa, Poland


\end{document}